\documentclass[a4paper,12pt]{article}
\usepackage[cp1250]{inputenc}
\usepackage[T1]{fontenc}
\usepackage{amsfonts,amssymb,amsmath,amsthm,bbding,bbm,enumitem,textcomp,url,yfonts,geometry}
\theoremstyle{plain}
\newtheorem{tw}{Theorem}[section]

\newtheorem{lemma}{Lemma}[section]
\newtheorem{cor}{Corollary}[section]

\theoremstyle{definition}
\newtheorem{definition}{Definition}[section]

\newtheorem{example}{Example}[section]
\newtheorem{cex}{Counterexample}[section]
\newtheorem{remark}{Remark}[section]

\theoremstyle{remark}

\newcommand{\semicopulas}{\mathfrak{S}}%class of semicopulas
\newcommand\by{\bar{y}}
\newcommand\bd{\bar{d}}

\newcommand{\bI}{\mathbf{I}}

\newcommand\cA{{\mathcal A}}

\newcommand\cF{{\mathcal F}} 
\newcommand\cM{{\mathcal M}}

\newcommand{\rA}{\mathrm{A}} 
\newcommand{\rB}{\mathrm{B}} 
\newcommand{\rS}{\mathrm{S}} %semicopul
\newcommand{\rM}{\mathrm{M}} %minimum
\newcommand{\rW}{\mathrm{W}} %semicopula Lukasiewicza

\newcommand{\cMm}[1][]{\cM^{#1}_{(X,\cA)}} % set of all monotone measure
\newcommand{\cFf}[1][\bar{y}]{\cF_{(X,\cA)}^{\,#1}} % set of measurable functions

\newcommand{\sP}{\mathsf{P}} %semicopul
\newcommand\vtr{\vartriangle}

\newcommand\st{\star}
\renewcommand\c{\circ}
\newcommand\lo{{\,\lozenge\,}} %zgadza siÄ™
\newcommand\loo{{\lozenge}}

\newcommand\calka[3][m]{\mathbf{I}_{#2}(#1,#3)}  %uogolniona seminormowa caĹ‚ka
\newcommand\calkab[3][m]{\mathbf{I}_{#2}\big(#1,#3\big)}  %uogolniona seminormowa caĹ‚ka z wiÄ™kszymi nawiasami

\newcommand\md{\,{\mathrm{d}}}

\renewcommand\ge{\geqslant}
\renewcommand\le{\leqslant}

\newcommand{\mI}[1]{\mathbbm{1}_{#1}}

\usepackage{xcolor}
\definecolor{darkgreen}{rgb}{0,0.5,0}

\title{General form of Chebyshev type inequality \\ for generalized Sugeno integral}
\author{Micha\l  \,\,Boczek, Anton Hovana, Ondrej Hutn\'{i}k
}
\date{ }
\begin{document}
\maketitle

\begin{abstract}
We prove a~general form of Chebyshev type inequality for generalized upper Sugeno integral in the form of necessary and sufficient condition. %This class of integrals includes  seminormed integrals (including the prominent Sugeno integral) as well as q-integrals. 
A key role in our considerations is played by the~class of $m$-positively dependent functions which includes comonotone functions as a~proper subclass. As a~consequence, we state an equivalent condition for Chebyshev type inequality to be true for all comonotone functions and any monotone measure. Our results generalize many others 
%MB bylo powtorzenie
obtained in the framework of q-integral, seminormed fuzzy integral and Sugeno integral on the real half-line.
Some further consequences of these results are obtained, among others Chebyshev type inequality for any functions. We also point out some flaws in existing results and provide their improvements.
\end{abstract}

\section{Introduction}
In many practical investigations, it is necessary to bound one quantity by another. The classical inequalities are very useful for this purpose. For instance, the classical Chebyshev integral inequality gives a lower bound for Lebesgue integral of product of two functions in terms of product of their Lebesgue integrals. The best result in the additive setting is as follows, see~\cite{Armstrong93}: \textit{If $(X,\cA)$ is a~measurable space,
then two $\cA$-measurable real functions $f$ and $g$ defined on $X$ satisfy the inequality
\begin{equation}\label{CH}
\int fg \md \sP \ge\int f\md \sP\int g\md \sP
\end{equation}
for any probability measure $\sP$ if and only if functions $f, g$ are comonotone.}

It is well-known that the classical integral inequalities (including the Chebyshev one) need not hold in general when replacing in~\eqref{CH} the probability measure by a non-additive measure and the additive (Lebesgue) integral by a non-additive integral. Nowadays, there is a~huge number of papers dealing with  inequalities for various non-additive integrals, see e.g. \cite{AgahiEslamiMohammadpourVaezpourYaghoobi12,AgahiMesiarOuyang09,KaluszkaOkolewskiBoczek14,OuyangMesiar09,YanOuyang19}. In this paper, we aim to provide a~general form of the Chebyshev type inequality (necessary as well as sufficient conditions) for the generalized upper Sugeno integral for $m$-positively dependent functions, as well as for any comonotone functions and monotone measure, see Section \ref{sec:Chebyshev}.  We show that our results generalize  many results from the literature, see Section \ref{sec:special}. 
%including the most recent results on Chebyshev type  inequality for q-integrals, see \cite{YanOuyang19}. 
Moreover, we improve some statements including the most recent results on Chebyshev type inequality for q-integral from \cite{YanOuyang19}  (see  Sections~\ref{sec:qint} and \ref{sec:seminormed}). Finally, in Section \ref{sub:further} we present Chebyshev type inequality  for \textit{all} functions under some mild assumptions on a~monotone measure.
We also present several examples demonstrating these results.

\section{Preliminaries}

Let $(X,\cA)$ be a~measurable space, where $\cA$ is a~$\sigma$-algebra of subsets of a~non-empty set $X.$ For a~given measurable space $(X,\cA)$ we denote the set of all $\cA$-measurable functions $f\colon X \to [0,\by]$ for some $\by\in (0,\infty]$ by $\cFf$. We also consider the set $\cMm$ of all \textit{monotone } (or, \textit{non-additive}) \textit{measures}, i.e., set functions
$m\colon \cA \to [0,\infty]$ satisfying $m(A)\le m(B)$ whenever $A\subset B$ with the boundary condition
$m(\emptyset) = 0$ and $m(X) > 0.$  If $m(X)=1,$ then $m$ is called a~\textit{capacity} and $\cMm[1]$ denotes the set of all capacities. Let $m(\cA\cap D)=\{m(A\cap D)\colon A\in \cA\}$ for a~fixed $D\in\cA$ and $m\in\cMm.$ To shorten the notation, in the case $D=X$ we denote by $m(\cA)$ the range of $m$.
%We say that $m$ is \textit{continuous} if $m\big(\bigcup_{n=1}^\infty A_n\big)=\lim\limits_{n\to \infty} m(A_n)$ and $m\big(\bigcap_{n=1}^\infty B_n\big)=\lim\limits_{n\to \infty} m(B_n)$ for $A_1\subset A_2\subset \ldots$ and $B_1 \supset B_2\supset \ldots$ with $A_n, B_n\in\cA$ for each $n$ and $m(B_k)<\infty$ for some $k.$

A binary operation $\circ\colon [0,\by]^2\to [0,\by]$ is called a~\textit{fusion function}. 
Pre-aggregation functions and aggregation functions are the most important examples of fusion functions (see \cite{BustinceBarrenecheaSesma-SaraLafuenteDimuroMesiarKolesarova17}). 
We say that a~function $\c\colon Y_1\times Y_2\to [0,\infty]$ is \textit{non-decreasing} if $a_1\c a_2\le b_1\c b_2$ whenever $a_i\le b_i,$ where $a_i,b_i\in Y_i\subset [0,\infty]$ for $i=1,2.$ 
%The binary operation $\c$ is \textit{increasing} if $a<b$ implies $a\c c<b\c c$ and $c\c a<c\c b.$
A~non-decreasing fusion function $\circ\colon [0,\by]^2\to [0,\by]$ is \textit{left-continuous} 
if it is left-continuous with respect to each coordinate.
With $\by=\infty$  the important fusion functions are $\rM(a,b) = a\wedge b$ and $\Pi(a,b)=ab$ (under the convention $0\cdot \infty=\infty\cdot 0=0$), where $a\wedge b=\min(a,b).$  Both of them are semicopulas for $\by=1$ as well. Recall that a~\textit{semicopula} (also called a~\textit{$t$-seminorm}) $\rS\colon [0,1]^2\to[0,1]$ is a~non-decreasing fusion function  such that $\rS(a,1)=\rS(1,a)=a$ (see \cite{Borzova-MolnarovaHalcinovaHutnik15,OuyangMesiar09}). 
From it follows that $\rS(x,0)= 0 = \rS(0,x)$ for all $x\in[0,1]$ and $\rS(x,y)\le x\wedge y$ for all $x,y\in[0,1]$. Another example is the~\textit{{\L}ukasiewicz semicopula} defined by $\rW(a,b)=(a+b-1)\vee 0,$ where $a\vee b=\max(a,b).$ 
%The  function $\rW$ is called the \textit{{\L}ukasiewicz semicopula}. 
The set of all semicopulas $\rS$ will be denoted by $\semicopulas.$

The \textit{generalized (upper) Sugeno integral} of $f\in \cFf$ on $D\in\cA$ with respect to $m\in\cMm$ is defined by 
\begin{align}
\label{defIcirc}
\calka{\c,D}{f} :=\sup_{t\in [0,\by]} \big\{t\c m(D\cap \{f\ge t\} )\big\},
\end{align}
where $\c\colon [0,\by] \times m(\cA\cap D)\to [0,\infty]$ is a~non-decreasing function  and $\{ f\ge t\} = \{x\in X\colon f(x)\ge t\}$ (see \cite{KaluszkaOkolewskiBoczek14}). To simplify the notation for $D=X$ we write 
\begin{align}\label{m01a}
\calka{\c}{f}:= \sup_{t\in[0, \by]} \big\{t\c m(\{f\ge t\})\big\}.
\end{align}
Note that 
$\calka{\c}{h\mI{D}}=\calka{\c,D}{h}$
for all $h\in\cFf$ and $m\in\cMm$ whenever $0\c m(X)=0\c m(D),$
where $\mI{A}$ denotes the indicator function of the~set $A.$ Replacing the operation $\c$ in the formula \eqref{m01a} with $\rM,$ $\Pi,$ $\rW$ and $\rS\in\semicopulas$ we get the \textit{Sugeno integral} $\bI_{\rM}$ \cite{SugenoPhD},  \textit{Shilkret integral} $\bI_\Pi$ \cite{Shilkret71}, \textit{opposite-Sugeno integral} $\bI_{\rW}$ \cite{Imanoka97} and \textit{seminormed fuzzy integral} $\bI_{\rS}$ \cite{Suarez-GarciaAlvarez-Gil86}, respectively. 

When investigating various inequalities for non-additive integrals, one usually has to restrict the class of measurable functions, or the class of monotone measures under consideration. In connection with the Chebyshev type inequalities, the class of $m$-positively dependent functions will play the key role. From now on, $h|_A$ denotes the~restriction of function $h$ to set $A.$

\begin{definition}\label{defmPD}
Let $f,g\in \cFf[k],$ $m\in\cMm$ and $A,B\in\cA$ with $k>0.$ Functions $f|_A$ and $g|_B$ are called \textit{$m$-positively dependent with respect to an operator $\vtr\colon m(\cA)\times m(\cA)\to m(\cA)$}, if
\begin{align}
\label{mu1}
m\big(A\cap B \cap \{ f\ge \alpha\}\cap \{ g\ge \beta\}\big)
\ge m(A\cap \{ f\ge \alpha\})
 \vtr m(B\cap \{ g\ge \beta\})
\end{align}
holds for all $\alpha,\beta\in [0,k].$
\end{definition}

The concept of $m$-positively dependent functions was introduced by Kaluszka et al. \cite{KaluszkaOkolewskiBoczek14} as a~natural generalization of positive dependency defined for the first time by Lehmann~\cite{Lehmann66} in probability theory. These ideas coincide for $\vtr=\cdot$ and probability measure $m=\sP.$  Also, comonotone functions on $A$ are included here being $m$-positively dependent with respect to $\vtr=\wedge$ with $B=A$ and any $m\in\cMm.$
Recall that two functions $f,g\colon X\to [0,\infty]$ are called \textit{comonotone on} $D\subset X$ if $(f(x)-f(y))(g(x)-g(y))\ge 0$  for all $x,y\in D.$
Now, we provide a~few examples of $m$-positively dependent functions.

\begin{example}\label{ex:new2}
All functions $f|_A,g|_B\in \cFf[\infty]$ are $m$-positively dependent with respect to any $\vartriangle$ whenever
$m\in \cMm$ satisfies the condition $m(C\cap D)\ge m(C)\vtr m(D)$ for all $C,D\in\cA.$  For instance, if $m\in\cMm$ is \textit{minitive}, i.e., $m(C\cap D)=m(C)\wedge m(D),$ with $m(X)\le 1,$ then any functions $f|_A,g|_B\in \cFf[\infty]$ are $m$-positively dependent for any $\vtr\in\semicopulas$ (e.g. $\rM$ or $\Pi$).
\end{example}

%\begin{example}
%Let $a\vtr b=b\mI{\{a+b>1\}}$ (the so-called G\"odel conjunction). Then all functions $f|_A$ and $g|_A$ are $m$-positively dependent for all $m\in\cMm$ and $A\in\cA$ such that $m(A)\le 0.5$.
%\end{example}

\begin{example}
Let $X=\{\omega_1,\omega_2\}$ and $m\in\cMm[1]$ be such that $m(\{\omega_2\})=p,$ where $p\in (0,1)$. 
Assume that $m(\{f\ge t\})=m(\{g\ge t\})=\mI{\{0\}}(t)+p\,\mI{(0,1]}(t)$ for all $t\ge 0.$
  Then $f|_X$ and $g|_X$ are $m$-positively dependent functions with respect to $\vtr$ such that $\vtr\le \cdot.$
%defined by $a\vtr b=a^{q_1}b^{q_2}$ for $q_1,q_2\ge 1.$
\end{example}

\begin{example}\label{ex:rec1}
Let $m\in\cMm[1]$ and $\vtr=\rW.$ 
%and $A,B\in\cA.$
Functions $f|_A,g|_B\in \cFf[\infty]$ are $m$-positively dependent if
\begin{align*}
m\big(A\cap B \cap \{ f\ge \alpha\}\cap \{ g\ge \beta\}\big)\ge \big(m(A\cap \{ f\ge \alpha\})+ m(B\cap \{ g\ge \beta\})-1\big)_+
\end{align*}
for any $\alpha,\beta\in [0,\infty],$
where $(a)_+=\max(a,0).$
The above inequality can be rewritten as follows
\begin{align*}
m^d\big((A\cap \{f\ge \alpha\})^c\cup (B\cap \{g\ge \beta\})^c\big)\le m^d\big((A\cap \{f\ge \alpha\})^c\big)+m^d\big((B\cap \{g\ge \beta\})^c\big)
\end{align*}
for all $\alpha,\beta,$ where $C^c=X\setminus C$ and $m^d(C)=1-m(C^c)$ is a~dual capacity. If $m^d$ is subadditive\footnote{A~monotone measure $m$ is \textit{subadditive}, if $m(A\cup B)\le m(A)+m(B)$ for all $A,B\in\cA.$}, then all functions $f|_A$ and $g|_B$ are $m$-positively dependent with respect to $\vtr$. For instance, all functions are $m$-positively dependent 
whenever
%if
$m$ is a~supermodular capacity\footnote{A~monotone measure $m$ is \textit{supermodular}, if $m(A\cup B)+m(A\cap B)\ge m(A)+m(B)$ for all $A,B\in\cA.$}. Note that a~distorted probability defined by $m(B)=h(\sP(B))$ for all $B\in\cA$ is supermodular. Here, $\sP$ is a~probability measure and $h\colon [0,1]\to [0,1]$ is an increasing and convex function such that $h(0)=0$ and $h(1)=1$, see \cite[p.\,17]{Denneberg94}.
%p.17
Distorted probabilities play the key role in  cumulative prospect theory \cite{KahnemanTversky79} and insurance \cite{KaluszkaKrzeszowiec12}.
\end{example}

\begin{example}\label{ex:godel}
Let $a\vtr b= a\otimes_G b:=b\mI{\{a>1-b\}}$ (the so-called G\"odel conjunction). Then all functions $f|_A,g|_A\in \cFf[\infty]$ are $m$-positively dependent for all $m\in\cMm[1]$  such that $m(\cA)\subset [0,0.5]$.
\end{example}

Further examples of $m$-positively dependent functions can be found in \cite{KaluszkaOkolewskiBoczek14}. When considering equality instead of inequality in \eqref{mu1} with $\vtr\in \mathfrak{S},$  $m\in\cMm[1]$ and $A=B=X,$ we may obtain a~connection between monotone measures and semicopulas. Under some quite common assumptions, one can derive an analogous result to Sklar's theorem (see \cite[Theorem 2.2.1]{DuranteSempi15}), i.e., every semicopula links the survival function\footnote{Function $\overline{F}_m\colon [0,1]^2\to [0,1]$ defined as $\overline{F}_m(x_1,x_2)=m((x_1,1]\cap (x_2,1])$ is called a~survival function associated with $m.$}  of the capacity to its marginal survival functions (see \cite[Theorem 8.3.3]{DuranteSempi15} and \cite[Theorem 9]{Scarsini96}).

\section{General Chebyshev type inequalities}\label{sec:Chebyshev}

Using the concept of $m$-positive dependency, we give a~new general form of Chebyshev type inequality for integral $\mathbf{I}_{\circ}$. The following lemma will be useful in the proof of our main result. Hereafter, we use the convention  $\sup_\emptyset=0.$ 

\begin{lemma}\label{lempom}
Let $\by\in (0,\infty]$ and $c, k\in(0,\by]$. Assume that $g\colon [0,c]\to [0,\by]$ and $h\colon [0,\by]\to [0,\by],$ where $g$ is a~non-decreasing function.
If $g$ is left-continuous and $\sup_{y\in [0,k]} h(y)\in [0,c]$, then $g(\sup_{y\in [0,k]} h(y))=\sup_{y\in [0,k]} g(h(y)).$
\end{lemma}
\proof
Clearly, $\sup_{y\in [0,k]} g(h(y))\le g(\sup_{y\in [0,k]} h(y)).$ On the other hand, let $(y_n)_1^\infty\in [0,k]$ be a~sequence such that  $h(y_n)\nearrow \sup_{y\in [0,k]} h(y),$ where $a_n\nearrow a$ means that a~non-decreasing sequence $(a_n)_1^\infty$ converges to $a.$  By  left-continuity and monotonicity of  $g$ we have
\begin{align*}
g(\lim_{n\to \infty} h(y_n))=\lim_{n\to\infty} g(h(y_n))\le \sup_{y\in [0,k]} g(h(y)).
\end{align*}
Hence, $g(\sup_{y\in [0,k]} h(y))=\sup_{y\in [0,k]} g(h(y)).$ \qed

\begin{tw}\label{ctw2}
Assume that $\by\in (0,\infty],$  $k\in (0,\by],$ $m\in\cMm,$ $\vtr\colon m(\cA)\times m(\cA)\to m(\cA)$ and $\st,\loo$ are non-decreasing fusion functions such that $\loo$ is left-continuous. Let $\varphi_i\colon [0,\by]\to [0,\by]$ be functions such that  $\varphi_1$ is non-decreasing and $\varphi_j$ are increasing and right-continuous, $\psi_i\colon [0,\varphi_i(\by)]\to [0,\by]$ be non-decreasing and $\psi_j$ be left-continuous,  where  $i=1,2,3$ and $j=2,3.$
Let $\c_i\colon [0,\by]\times m(\cA)\to [0,\by]$ be non-decreasing such that $\varphi_i(\by)\c_i m(X)\le\varphi_i(\by)$ and
$\by\c_j 0=0$
%$\varphi_j(\by)\c_j 0=0$
with $i=1,2,3$ and $j=2,3.$  

(a) Assume that
\begin{align}\label{c2a}
\psi_1(\varphi_1(a\st b) \c_1 (c\vtr d))\ge \psi_2(\varphi_2(a)\c_2 c)\lo \psi_3 (\varphi_3(b)\c_3 d)
\end{align}
holds for all $a,b\in [0,k]$ and $c,d\in m(\cA)$. If  $f|_A,g|_B\in\cFf[k]$
are $m$-positively dependent with respect to $\vtr$, then
\begin{align}\label{c2b}
\psi_1\big(\calka{\c_1,A\cap B}{\varphi_1(f\st g)}\big)\ge \psi_2\big(\calka{\c_2,A}{\varphi_2(f)}\big)\lo \psi_3\big(\calka{\c_3,B}{\varphi_3(g)}\big).
\end{align}

(b) Suppose that $\by\c_1 0=0,$ $a\vtr 0=0\vtr a=0$ for all $a\in m(\cA)$ and the condition $(Z_1)$ holds, i.e., for all $c,d\in m(\cA)$ there exists sets $C,D\in\cA$ such that $c=m(C)$, $d=m(D)$ and $m(C\cap D)=m(C)\vtr m(D).$
If \eqref{c2b} is true for all $m$-positively dependent functions $f|_A,g|_B\in\cFf[k]$ with respect to $\vtr$,  then \eqref{c2a} holds for all $a,b\in [0,k]$ and $c,d\in m(\cA).$
\end{tw}

\proof
(a) Firstly, note that $\psi$'s are well-defined
%MB: previously well defined
in \eqref{c2a} as $\varphi_j(a)\c_j c\in [0,\varphi_j(\by)]$ for all $a\in [0,k]$ and $c\in m(\cA)$, where $j=2,3$, and $\varphi_1(a\st b)\c_1 (c\vtr d)\in  [0,\varphi_1(\by)]$ for all $a,b\in [0,k]$ and $c,d\in m(\cA)$, since $\varphi_i(\by)\c_i m(X)\le \varphi_i(\by)$ for $i=1,2,3.$ Then $m$-positive dependency of $f|_A$ and~$g|_B$ implies that
\begin{align*}
m(A\cap B \cap \{ f\ge a\}\cap \{ g\ge b\})
\ge m(A\cap \{ f\ge a\})\vtr m(B\cap \{ g\ge b\})
\end{align*}
for all $a,b\in [0,k].$
The monotonicity of $\st$ and  $m$ yields
\begin{align*}
m(A\cap B \cap \{ f\st g\ge a\st b\}) \ge m(A\cap \{ f\ge a\})\vtr m( B \cap\{ g\ge b\})
\end{align*}
for any $a,b\in [0,k].$ As $\varphi_1$ is non-decreasing and $\varphi_j$ are increasing for $j=2,3$, we get
\begin{align*}
m\big(A\cap B \cap \lbrace \varphi_1(f\st g)
\ge \varphi_1(a\st b)\rbrace\big)\ge c_a\vtr d_b
\end{align*}
for all $a,b\in [0,k]$  
with $c_a:=m(A\cap \lbrace \varphi_2(f)\ge \varphi_2(a)\rbrace)$ and $d_b:=m(B\cap \{ \varphi_3(g)\ge \varphi_3(b)\}).$ 
From the assumption on monotonicity of $\c_1$ and $\psi_1$, we obtain
\begin{align*}
\psi_1\big(\varphi_1(a\st b)\c_1 m(A\cap B\cap \{\varphi_1(f\st g)\ge \varphi_1(a\st b)\})\big)
\ge \psi_1(\varphi_1(a\st b)\c_1 (c_a \vtr d_b))
\end{align*}
for any $a,b\in[0,k]$.
By the definition of integral $\mathbf{I}_{\c_1},$ function $\psi_1$ and from \eqref{c2a}, we get
\begin{align}\label{c4b}
\psi_1\big(\calka{\c_1, A\cap B}{\varphi_1(f\st g)}\big)\ge\psi_2(\varphi_2(a)\c_2 c_a)\lo \psi_3(\varphi_3(b)\c_3d_b)
\end{align}
for any $a,b\in [0,k].$
For a~fixed $y\in [0,\by]$ put
$H(x)=x\lo y$  for all $x\in [0,\by]$. Also, put $h(a)=\varphi_2(a)\c_2 c_a$ for all $a\in [0,k].$
The function $H(\psi_2)\colon [0,\varphi_2(\by)]\to [0,\by]$ is non-decreasing and left-continuous as $\loo$ and $\psi_2$ are non-decreasing and left-continuous. Then Lemma~\ref{lempom} yields
\begin{align}\label{c5a}
\sup_{a\in [0,k]} H\big(\psi_2(h(a))\big)= H\big(\psi_2(\sup_{a\in [0,k]} h(a))\big).
\end{align}
Moreover, $\sup_{a\in [0,k]} h(a)=\sup_{a\in [0,\by]} h(a)$ due to the fact that $\varphi_2(\by)\c_2 0=0$ and $f\le k\le \by.$  Putting $t=\varphi_2(a),$ we get
\begin{align*}
\sup_{a\in [0,k]} h(a) & = \sup_{t \in \varphi_2([0,\by])}\big\{ t\c_2 m(A\cap\{\varphi_2(f)\ge t\})\big\},
\end{align*}
where $\varphi_2 ([0,\by])$ is the image of $\varphi_2.$
Note that
\begin{align*}
\sup_{a\in [0,k]} h(a)=\sup_{t \in [0,\varphi_2(0))}\{ t\c_2 m(A)\}\vee \sup_{t \in \varphi_2([0,\by])}\big\{ t\c_2 m(A\cap\{\varphi_2(f)\ge t\})\big\}\vee \sup_{t \in (\varphi_2(\by),\by]} \{ t\c_2 0\}
\end{align*}
since $\sup_{t \in [0,\varphi_2(0))}\{ t\c_2 m(A)\}=\varphi_2(0)\c_2 m(A)$ and $\by\c_2 0=0.$
From right-continuity of $\varphi_2$  and monotonicity of $\c_2$
we have
\begin{align}\label{c5e}
\sup_{a\in [0,k]} h(a) =\sup_{t \in [0,\by]} \big\{t\c_2 m(A\cap\{\varphi_2(f)\ge t\})\big\}.
\end{align}
From \eqref{c4b}-\eqref{c5e} we obtain
\begin{align*}
\psi_1\big(\calka{\c_1, A\cap B}{\varphi_1(f\st g)}\big)
\ge \psi_2\big(\calka{\c_2,A}{\varphi_2(f)}\big) \lo \psi_3(\varphi_3(b)\c_3 d_b)
\end{align*}
for all $b\in [0,k].$ Proceeding similarly with the supremum in $b\in [0,k]$, we get \eqref{c2b}.

(b) Fix $a,b\in [0,k]$ and $c,d\in m(\cA).$  Define $f=a\mI{X}$ and $g=b\mI{X}$, and consider $A,B\in \cA$ satisfying the condition $(Z_1)$.
Then $f|_A$ and $g|_B$ are $m$-positively dependent with respect to $\vtr$. Indeed,
\eqref{mu1} takes the form
$0=m(\emptyset)\ge m(C\cap D) = m(C)\vtr m(D)=0$
for $\alpha>a$ or $\beta>b,$ where $C=\{f|_A\ge \alpha\}$ and $D=\{g|_B\ge \beta\},$ as $0\vtr x=x\vtr 0=0.$  If $\alpha\le a$ and $\beta\le b,$ then we have $m(A\cap B) = m(A)\vtr m(B),$ so the inequality \eqref{mu1} is satisfied.
Thus,
\begin{align*}
\calka{\c_2,A}{\varphi_2(f)}& = \varphi_2(a) \c_2 m(A), \\
\calka{\c_3,B}{\varphi_3(g)} & = \varphi_3(b)\c_3 m(B),\\
\calka{\c_1,A\cap B}{\varphi_1(f\st g))} & = \varphi_1(a\st b)\c_1 m(A\cap B)
=\varphi_1(a\st b)\c_1 (m(A)\vtr m(B)),
\end{align*}
since $\by\c_i 0=0$ for all $i.$ Applying \eqref{c2b} we obtain \eqref{c2a} for all $a,b\in [0,k]$ and $c,d\in m(\cA).$
\qed\medskip

%MB: nie usuwac, moze przydac sie kiedys tam
%Aby udowodnić równanie \eqref{c5e} można zastapic warunek $y\c_j 0=0$ na warunek $\varphi_j(\by)=\by$ o ile $k=\by.$

Note that Theorem~\ref{ctw2} generalizes Theorem 2.3 from \cite{KaluszkaOkolewskiBoczek14}. Moreover, the inequality \eqref{c2a} has been investigated for some functions in~\cite{KaluszkaOkolewskiBoczek14}.  
Now, we mention two examples where the inequality \eqref{c2b} becomes equality.
\begin{example}\label{ex:rown1}
The equality in \eqref{c2b} holds if $k=\by=1,$ $\c_i=\wedge$, $\st=\loo=\cdot$,  $\psi_i(x)=x^{1/p_i},$ $\varphi_i(x)=x^{p_i},$ $f=a\mI{D}$ and $g=b\mI{D},$ where $a,b\in [0,1],$ $p_i>0,$ $D\subset A\cap B$ and $(ab)^{p_1}\vee a^{p_2}\vee b^{p_3}\le m(D)$ for $i=1,2,3.$
\end{example}

\begin{example}\label{ex:rown2}If $k=\by=1,$ $\c_i\in\semicopulas$, $A=B=X,$ $\st=\loo=\cdot$, $f=a\mI{X},$  $g=b\mI{X},$ $\psi_i(x)=x^{q}$ and $\varphi_i(x)=x^p,$ where $a,b\in [0,1],$ $p,q>0$ and $i=1,2,3,$  then the equality in \eqref{c2b} holds for $m\in\cMm[1].$ 
\end{example}

The condition $(Z_1)$ appears naturally in the context of $m$-positively dependent functions. In order to avoid  it,
%this condition,
we have to consider the class of comonotone functions on the whole $X$ as a~special case of $m$-positively dependent functions with respect to $\vtr=\wedge.$ Firstly, we give the following useful lemma.

\begin{lemma}\label{lem:M2}
Let $\by\in (0,\infty]$, $k\in (0,\by]$ and $\st,\loo$ be  fusion functions such that $\loo$ is non-decreasing. 
Assume that $\varphi_i\colon [0,\by]\to [0,\by]$  and $\psi_i\colon [0,\varphi_i(\by)]\to [0,\by]$ are non-decreasing for $i=1,2,3$. Let $\c_i\colon [0,\by]\times D\to [0,\by]$ be non-decreasing such that
$\varphi_i(\by)\c_i \bd\le \varphi_i(\by)$ for $i=1,2,3,$ where $D\subset [0,\infty]$ and $\bd=\sup D\in D.$
Then the following conditions are equivalent:
\begin{itemize}[noitemsep]
\item[$(C_1)$]   $\psi_1(\varphi_1(a\st b)\c_1(c\wedge d))\ge \psi_2(\varphi_2(a)\c_2 c)\lo \psi_3 (\varphi_3(b)\c_3 d)$ for all $a,b\in [0,k]$ and $c,d\in D;$
\item[$(C_2)$] $\psi_1(\varphi_1(a\st b)\c_1 c)\ge \big[\psi_2(\varphi_2(a)\c_2 c)\lo \psi_3(\varphi_3(b)\c_3 \bd)\big]\vee \big[\psi_2(\varphi_2(a)\c_2 \bd)\lo \psi_3(\varphi_3(b)\c_3 c)\big]$  for all $a,b\in[0,k]$ and $c\in D$.
\end{itemize}
\end{lemma}

\proof Similar as in the proof of Theorem \ref{ctw2}, one can show that $\psi$'s are well-defined in $(C_1)$ and $(C_2),$ as $\varphi_i(\by)\c_i \bd\le \varphi_i(\by)$ for $i=1,2,3.$

``$(C_1) \Rightarrow (C_2)$''
Putting $d=\bd$ in $(C_1)$ and then $c=\bd$ in $(C_1)$ we get $(C_2)$.\\
\noindent ``$(C_2) \Rightarrow (C_1)$'' By $(C_2)$  and by monotonicity of $\c_3,$ $\psi_3$ and $\loo,$ we get
\begin{align}\label{m1}
\psi_1(\varphi_1(a\st b)\c_1 c) \ge  \psi_2(\varphi_2(a)\c_2 c)\lo \psi_3(\varphi_3(b)\c_3\bd)\ge  \psi_2(\varphi_2(a)\c_2 c)\lo\psi_3(\varphi_3(b)\c_3 d)
\end{align} for all $a,b\in [0,k]$ and $c,d\in D.$ 
Similarly, we obtain the following inequalities
\begin{align}\label{m2}
\psi_1(\varphi_1(a\st b)\c_1 d)\ge  \psi_2(\varphi_2(a)\c_2 \bd)\lo \psi_3(\varphi_3(b)\c_3 d)\ge  \psi_2(\varphi_2(a)\c_2 c)\lo\psi_3(\varphi_3(b)\c_3 d)
\end{align} for any $a,b\in [0,k]$ and $c,d\in D.$ Combining \eqref{m1} and \eqref{m2} with monotonicity of $\c_1$ and $\psi_1$ we get the condition $(C_1)$.
\qed\medskip

Now, we can state the second version of Chebyshev type inequality omitting the condition $(Z_1).$

\begin{tw}\label{ctw3}
Assume that $\by\in (0,\infty],$ $k\in (0,\by]$ and $\st,\loo$ are non-decreasing fusion functions such that $\loo$ is left-continuous. Let $\varphi_i\colon [0,\by]\to [0,\by]$ be functions such that  $\varphi_1$ is non-decreasing and $\varphi_j$ are increasing and right-continuous, $\psi_i\colon [0,\varphi_i(\by)]\to [0,\by]$ be non-decreasing and $\psi_j$ be left-continuous, where  $i=1,2,3$ and $j=2,3.$
Let $\c_i\colon [0,\by]^2\to [0,\by]$ be non-decreasing such that $\varphi_i(\by)\c_i \by\le\varphi_i(\by)$ and
$\by\c_j 0=0,$
%$\varphi_j(\by)\c_j 0=0,$
where $i=1,2,3$ and $j=2,3.$

(a) Suppose that the inequality
\begin{align}\label{m3a}
\psi_1(\varphi_1(a\st b) \c_1 (c\wedge d))\ge \psi_2(\varphi_2(a)\c_2 c)\lo \psi_3 (\varphi_3(b)\c_3 d)
\end{align}
holds for all $a,b\in [0,k]$ and $c,d\in [0,\by].$ If $f,g\in\cFf[k]$ are comonotone on $X$ and $m\in\cMm$ such that $m(X)\le \by,$ 
then
\begin{align}\label{m3b}
\psi_1\big(\calka{\c_1}{\varphi_1(f\st g)}\big)
\ge \psi_2\big(\calka{\c_2}{\varphi_2(f)}\big)\lo \psi_3\big(\calka{\c_3}{\varphi_3(g)}\big).
\end{align}

(b)  Assume that $0\st k=k\st 0=0,$ $\by\c_1 0=0$ and $\varphi_i(0)\c_i \by=0$ for $i=1,2,3.$  If \eqref{m3b} is true for all comonotone functions $f,g\in\cFf[k]$ on $X$ and any $m\in\cMm$ such that $m(X)\le \by,$ then \eqref{m3a} holds for all $a,b\in [0,k]$ and $c,d\in [0,\by].$
\end{tw}

\proof
Part (a) can be proved in the same way as  Theorem~\ref{ctw2}\,(a) with $\vtr=\wedge$ and $A=B=X,$ so we omit it.

(b)  Define $f=a\mI{A}$ and $g=b\mI{X},$ where $a,b\in [0,k]$ and $A\in\cA.$ Then
\begin{align*}
\calka{\c_2}{\varphi_2(f)}& =(\varphi_2(0)\c_2 m(X))\vee ( \varphi_2(a) \c_2 m(A))\vee (\by \c_2 0)
=\varphi_2(a) \c_2 m(A),
\end{align*}
 since $\varphi_2(0)\c_2 m(X) =0$ and $\by \c_2 0=0.$
 Similarly, one can check that
\begin{align*}
\calka{\c_3}{\varphi_3(g)} & = \varphi_3(b)\c_3 m(X),\\ \calka{\c_1}{\varphi_1(f\st g)} & =\varphi_1(a\st b)\c_1 m(A),
\end{align*}
as $0\st k=0.$
Applying \eqref{m3b} we obtain 
\begin{align*}
\psi_1(\varphi_1(a\st b)\c_1 m(A))
\ge \psi_2(\varphi_2(a)\c_2 m(A))\lo \psi_3(\varphi_3(b)\c_3 m(X)).
\end{align*}
Similarly, for $f=a\mI{X}$ and $g=b\mI{A}$ we get the inequality
\begin{align*}
\psi_1(\varphi_1(a\st b)\c_1 m(A))\ge \psi_2(\varphi_2(a)\c_2 m(X))\lo \psi_3(\varphi_3(b)\c_3 m(A)).
\end{align*}
In consequence, from \eqref{m3b} for any comonotone functions and any $m$ such that $m(X)\le \by,$ we get
\begin{align*}
\psi_1(\varphi_1(a\st b)\c_1 c)
\ge \big[\psi_2(\varphi_2(a)\c_2 c)\lo \psi_3(\varphi_3(b)\c_3 \by)\big]
\vee \big[\psi_2(\varphi_2(a)\c_2 \by)\lo \psi_3(\varphi_3(b)\c_3 c)\big]
\end{align*}  for all $a,b\in[0,k]$ and $c\in [0,\by].$
Using Lemma~\ref{lem:M2} with $D=[0,\by]$ we obtain \eqref{m3a}.
\qed\medskip

The next example demonstrates that in some cases we cannot use Theorem~\ref{ctw3}\,(a), but Theorem~\ref{ctw2}\,(a) still works.

\begin{example} 
Put $k=\by=1,$ $A=B=X,$  $\varphi_i(x)=\psi_i(x)=x$ for all $x\in[0,1],$ $\loo=\st=\cdot,$ $\c_i=\rW,$ $m(\cA)=\{0,1\}$ and $\vtr=\wedge$ in Theorem \ref{ctw2}\,(a), where $i=1,2,3.$ Then the inequality \eqref{c2a} takes the form
$$\rW(ab,c\wedge d)\ge \rW(a,c)\,\rW(b,d)$$ for all $a,b\in[0,1]$ and $c,d\in \{0,1\}.$ From Theorem~\ref{ctw2}\,(a), the Chebyshev type inequality
\begin{align*}
\calka{\rW}{fg}\ge \calka{\rW}{f}\, \calka{\rW}{g}
\end{align*}
holds for any comonotone functions $f,g\in\cFf[1]$ on $X.$
However, we cannot use Theorem~\ref{ctw3}\,(a), since  
$$\rW(ab,c\wedge d)\ge \rW(a,c)\,\rW(b,d)$$ is not valid
for all $a,b,c,d\in[0,1].$ Indeed, it is enough to take $a=b=0{.}5$ and $c=d=0{.}75.$
\end{example}

\section{Special cases}\label{sec:special}

In this section, we derive some consequences of theorems proven in the previous section for several classes of integrals which are known in the literature.  

\subsection{Q-integral}\label{sec:qint}
Recall that a~non-decreasing fusion function $\otimes \colon [0,1]^2\to [0,1]$ is said to be a~\textit{fuzzy conjunction}
if $1\otimes 1=1$ and $0\otimes  1=1\otimes  0=0\otimes 0=0$  (see \cite{DuboisPradeRicoTeheux17}). The most important examples of fuzzy conjunction are: 
the G\"odel conjunction $\otimes_G$ (see Example~\ref{ex:godel})
%$a\otimes_{G} b := b\mI{\{a>1-b\}}(a,b)$ 
and the contrapositive G\"odel conjunction $a\otimes_{GC} b = a\mI{\{a>1-b\}}(a,b).$ Dubois et al. \cite{DuboisPradeRicoTeheux17} introduced and studied the $q$-integral defined as
\begin{align}
\int_m^\otimes f=\sup_{t\in [0,1]}\big\{m (\{f\ge t\})\otimes t\big\},
\end{align} 
where $\otimes$ denotes a~fuzzy conjunction and $(m,f)\in\cMm[1]\times\cFf[1].$ This definition is motivated by alternative ways of using weights of qualitative criteria in min- and max-based aggregations, that make intuitive sense as tolerance thresholds. Note that the research on Chebyshev type inequality for q-integral has been initiated by Kaluszka et al. \cite{KaluszkaOkolewskiBoczek14} even before its formal definition by Dubois \cite{DuboisPradeRicoTeheux17}, see for example \cite[Theorem 2.1]{KaluszkaOkolewskiBoczek14} for $Y=[0,1],$ $\mu\in\cMm[1]$ and $a\c_i b=b\otimes a.$ From this point of view, the claim in~\cite{YanOuyang19} about starting the research by the authors is misleading.

\begin{cor}\label{cor_ouyang}
 Let $\otimes\colon [0,1]^2\to [0,1]$ be a~fuzzy conjunction and $\st\colon [0,1]^2\to [0,1]$ be non-decreasing, left-continuous and $1\st 0=0=0\st 1.$ Suppose that $\varphi_i\colon [0,1]\to [0,1]$ is continuous and increasing such that $\varphi_i(0)=0$  and $1\otimes \varphi_i(1)\le \varphi_i(1)$ for $i=1,2,3.$ Then the following statements are equivalent:
\begin{itemize}[leftmargin=0.6cm]
\item[(i)]  $\varphi_1^{-1}(a\otimes  \varphi_1(b \st c))
\ge  \big[\varphi_2^{-1}(a \otimes \varphi_2(b)) \st \varphi_3^{-1}(1 \otimes \varphi_3(c))\big] \vee\big[\varphi_2^{-1}(1 \otimes \varphi_2(b)) \st \varphi_3^{-1}(a \otimes \varphi_3(c))\big]$ for any $a,b,c \in [0,1];$
\item[(ii)] $\varphi^{-1}_1\big(\int_m^\otimes \varphi_1(f\st g)\big)\ge \varphi_2^{-1}\big(\int_m^\otimes \varphi_2(f)\big)\st \varphi_3^{-1}\big(\int_m^\otimes \varphi_3(g)\big)$ holds for any comonotone functions $f,g\in \cFf[1]$ on $X$ and any capacity $m.$
\end{itemize}
\end{cor}
\proof
Put $\loo=\st,$ $\by=k=1,$ $a\c_i b=b\otimes a$ and $\psi_i=\varphi_i^{-1}$ for $i=1,2,3$ in Theorem~\ref{ctw3} and Lemma~\ref{lem:M2}. Functions $\psi_i$ are well-defined as $\varphi_i(0)=0.$ Using Lemma \ref{lem:M2} with $D=[0,1]$ and Theorem~\ref{ctw3}, we get the statement.
\qed\medskip

Note that if we put $\varphi_i(x)=x$ in Corollary~\ref{cor_ouyang}, we get \cite[Theorem 3.5]{YanOuyang19}. 

\begin{example}
Put $\varphi_i=\varphi$ and $\st=\cdot$ in Corollary~\ref{cor_ouyang} where $\varphi\colon [0,1]\to [0,1]$ is continuous and increasing such that $\varphi(0)=0$ and $\varphi(1)=1.$ We show that the following inequality 
\begin{align}\label{b2}
\varphi^{-1}(a\otimes  \varphi(b c))\ge  \big[\varphi^{-1}(a \otimes \varphi(b))\cdot  \varphi^{-1}(1 \otimes \varphi(c))\big]\vee\big[\varphi^{-1}(1 \otimes \varphi(b)) \cdot \varphi^{-1}(a \otimes \varphi(c))\big]  
\end{align}
does not hold for any $a,b,c\in [0,1]$ and $\otimes\in \{\otimes_{G},\otimes_{GC}\}.$
Putting $b=1$ in \eqref{b2}, we get
\begin{align*}
\varphi^{-1}(a\otimes  \varphi(c))\ge  \big[\varphi^{-1}(a \otimes 1) \cdot \varphi^{-1}(1 \otimes \varphi(c))\big] \vee\varphi^{-1}(a \otimes \varphi(c)) 
\end{align*}
for all $a,c\in [0,1].$ Consider $a>0$ and $\varphi(c)>0$ such that $a+\varphi(c)\le 1.$ Then we have
\begin{align*}
0\ge  \varphi^{-1}(a \otimes 1) \cdot \varphi^{-1}(1 \otimes \varphi(c))
\end{align*}
for $\otimes\in \{\otimes_{G},\otimes_{GC}\}$ which leads to the contradiction.  In consequence, the Chebyshev type inequality
\small{$$\varphi^{-1}\Bigg(\int_m^\otimes \varphi(f\st g)\Bigg)\ge \varphi^{-1}\Bigg(\int_m^\otimes \varphi(f)\Bigg)\st \varphi^{-1}\Bigg(\int_m^\otimes \varphi(g)\Bigg)$$}\normalsize
does not hold for any comonotone functions $f,g\in\cFf[1]$ on $X$ and any capacity $m$ for $\otimes\in \{\otimes_{G},\otimes_{GC}\}.$
\end{example}

\subsection{Seminormed fuzzy integral}\label{sec:seminormed}
Recall that for $\circ = \rS\in\mathfrak{S}$ (a~semicopula) the generalized Sugeno integral $\mathbf{I}_{\circ}$ coincides with the seminormed fuzzy integral being the smallest semicopula-based universal integral (see~\cite{Borzova-MolnarovaHalcinovaHutnik15}). 
Observe that the seminormed fuzzy integral is a~special case of q-integral with the conjunction replaced by the semicopula in the following way $a\otimes b= \rS(b,a).$ Applying Theorem~\ref{ctw2} with $\by=1$ and $\c_i=\rS_i$ for $i=1,2,3$, we get the following Chebyshev type inequality for seminormed fuzzy integral.

\begin{cor}
Assume that $m\in\cMm[1],$ $\vtr\colon m(\cA)\times m(\cA)\to m(\cA),$ $k\in (0,1],$  $\rS_i\in\semicopulas$ for $i=1,2,3$ and $\st,\loo\colon [0,1]^2\to [0,1]$ are non-decreasing and $\loo$ is left-continuous. Let $\varphi_i\colon [0,1]\to [0,1]$ be functions such that  $\varphi_1$ is non-decreasing and $\varphi_j$ are increasing and right-continuous, $\psi_i\colon [0,\varphi_i(1)]\to [0,1]$ be non-decreasing  and $\psi_j$ be left-continuous, where  $i=1,2,3$ and $j=2,3.$

(a) Suppose that
\begin{align}\label{d2a}
\psi_1\big(\rS_1(\varphi_1(a\st b),c\vtr d)\big)\ge \psi_2\big(\rS_2(\varphi_2(a),c)\big)\lo \psi_3 \big(\rS_3(\varphi_3(b),d)\big)
\end{align}
holds for all $a,b\in [0,k]$ and $c,d\in m(\cA)$. If  $f|_A,g|_B\in\cFf[k]$
are $m$-positively dependent with respect to $\vtr,$ then
\begin{align}\label{d2b}
\psi_1\big(\calka{\rS_1,A\cap B}{\varphi_1(f\st g)}\big)\ge \psi_2\big(\calka{\rS_2,A}{\varphi_2(f)}\big)\lo \psi_3\big(\calka{\rS_3,B}{\varphi_3(g)}\big).
\end{align}

(b) Assume that $a\vtr 0=0\vtr a=0$ for all $a\in m(\cA)$ and the condition $(Z_1)$ holds.
If \eqref{d2b} is true for all $m$-positively dependent functions $f|_A,g|_B\in \cFf[k]$ with respect to $\vtr,$ then \eqref{d2a} holds for any $a,b\in [0,k]$ and $c,d\in m(\cA).$
\end{cor}

Next corollary gives a~necessary and sufficient condition for Chebyshev type inequality for $\mathbf{I}_\rS$ for comonotone functions on $X$ and any monotone measure.

\begin{cor}\label{corAEMV12}
Assume $\rS\in\semicopulas$ and operation $\st\colon [0,1]^2\to [0,1]$ is  non-decreasing and  left-continuous such that $0\st 1=1\st 0=0.$ Let  $\varphi_i\colon [0,1]\to [0,1]$ be functions such that  $\varphi_1$ is non-decreasing and $\varphi_j$ are increasing and right-continuous, $\psi_i\colon [0,\varphi_i(1)]\to [0,1]$ be non-decreasing  and $\psi_j$ be left-continuous, where  $i=1,2,3$ and $j=2,3.$ Then the Chebyshev type inequality
\begin{align}\label{b4a}
\psi_1\big(\calkab{\rS}{\varphi_1(f\st g)}\big)\ge \psi_2\big(\calka{\rS}{\varphi_2(f)}\big)\st \psi_3\big(\calka{\rS}{\varphi_3(g)}\big)
\end{align}
is fulfilled for any comonotone functions $f,g\in\cFf[1]$ on $X$ and any $m\in\cMm[1]$ if and only if
\begin{align}\label{b4b}
\psi_1\big(\rS(\varphi_1(a\st b),c)\big)\ge
\big[\psi_2\big(\rS(\varphi_2(a),c)\big)\st b\big]\vee \big[a\st \psi_3\big(\rS(\varphi_3(b),c)\big)\big]
\end{align}
holds for all $a,b,c\in [0,1].$
\end{cor}
\proof
Use Lemma~\ref{lem:M2} with $k=\by=1,$ $\loo=\st,$ $D=[0,1]$  and $\c_i=\rS,$ and Theorem~\ref{ctw3}. 
\qed\medskip

Corollary~\ref{corAEMV12} with  $\psi_i=\varphi_i^{-1}$  is an improvement of \cite[Theorem 4.1 with $A=X$]{AgahiEslamiMohammadpourVaezpourYaghoobi12}, where the integral inequality \eqref{b4a} has been proved whenever the condition \eqref{b4b} holds. Setting $\psi_i(x)=\varphi_i(x)=x$ in Corollary~\ref{corAEMV12}, we get the following result.

\begin{cor}\label{cor:qint}
Let $\rS\in\semicopulas,$ $\st\colon [0,1]^2\to [0,1]$  be non-decreasing and left-continuous such that $0\st 1=1\st 0=0.$
The Chebyshev type inequality 
\begin{align}\label{b6}
\calkab{\rS}{f\st g}\ge \calka{\rS}{f}\st \calka{\rS}{g}
\end{align}
is satisfied for any comonotone functions $f,g\in\cFf[1]$ on $X$ and any $m\in\cMm[1]$ if and only if
\begin{align}\label{b7}
\rS(a\st b,c)\ge
(\rS(a,c)\st b)\vee (a\st \rS(b,c))
\end{align}
is valid for all $a,b,c\in [0,1].$
\end{cor}

Corollary~\ref{cor:qint} improves \cite[Theorem 3.1 with $A=X$]{OuyangMesiar09}.
In \cite[Corollary 3.6]{YanOuyang19} authors claim (without any proof) that the inequality \eqref{b6} is true for any comonotone functions $f,g\in\cFf[1]$ on $X$ and any capacity $m$ if and only if
\begin{align}\label{b8}
 \rS(c,a\st b)\ge
(\rS(c,a)\st b)\vee (a\st \rS(c,b))
\end{align}
for any $a,b,c\in [0,1].$
Clearly, condition \eqref{b7} coincides with \eqref{b8} only if $\rS\in \semicopulas$ is commutative.

Furthermore, in~\cite{DarabyGhadimi14} we can find an incorrect statement about Chebyshev type inequality with comonotone functions, since the sufficient condition used therein does not have the form \eqref{b4b}. Below, we present a~counterexample to \cite[Theorem 3.5]{DarabyGhadimi14}. Here we use the same notation of functions as in the original paper~\cite{DarabyGhadimi14}.

\begin{cex}\label{cex2}
Put $\mathrm{T}=\st=\rW$, $f=0{.}5\,\mI{A}$, $g=0{.}8\,\mI{A}$, $\mu(A)=0{.}9$ and $s=2$ in \cite[Theorem 3.5]{DarabyGhadimi14}. Clearly, functions $f$ and $g$ are comonotone. Then the sufficient condition (3.4) from~\cite{DarabyGhadimi14} in the form
\begin{align*}
\big((a+b-1)_++c-1\big)_+\ge \big((a+c-1)_++b-1\big)_+\vee \big(a+(b+c-1)_+-1\big)_+
\end{align*}
is satisfied for all $a,b,c\in [0,1].$ After a~simple calculation we get
\begin{align*}
\sqrt{\calka[\mu]{\rW,A}{\rW^2(f,g)}}&=\sqrt{\rW(0{.}3^2,0{.}9)}=0,\\
\rW\Big(\sqrt{\calka[\mu]{\rW,A}{f^2}}, \sqrt{\calka[\mu]{\rW,A}{g^2}}\Big)&=\Big(\sqrt{\rW(0{.}5^2,0{.}9)}+\sqrt{\rW(0{.}8^2,0{.}9)}-1\Big)_+\approx 0{.}122145,
\end{align*}
which contradicts the Chebyshev type inequality (3.5) from Theorem 3.5 in \cite{DarabyGhadimi14}.
\end{cex}

\subsection{Sugeno integral}\label{sec:sug}

Chebyshev type inequalities for Sugeno integral $\mathbf{I}_{\rM}$  in case of functions from $\cFf[1]$ can be obtained from the corresponding results for $\mathbf{I}_{\rS}$ (with the Sugeno integral regarding as a~special case of seminormed integral), see e.g. Corollary~\ref{corAEMV12}. On the other hand, Chebyshev type inequalities for Sugeno integral in case of functions from $\cFf[\by],$ $\by>0,$ can be obtained using Theorem~\ref{ctw2} or Theorem~\ref{ctw3} in a~similar manner as in the previous subsection. However, this may lead to some functional inequalities that need not be easy to verify.
In the case of Sugeno integral and comonotone functions, the functional inequality may be replaced by easier condition $\st\le \wedge.$ We start with the helpful lemma.

\begin{lemma}\label{lempom2}
Let $D\subset [0,\infty]$ and $\st$ be a~non-decreasing fusion function.
Assume that $\varphi_i\colon [0,\by]\to [0,\by]$   and $\psi_i\colon [0,\varphi_i(\by)]\to [0,\by]$ are non-decreasing  for $i=1,2,3$ such that  $\varphi_1(\by)= \varphi_j(\by),$ $\psi_1\ge \psi_j$ and $\psi_j(\varphi_j(x))\le  x\le \psi_1(\varphi_1(x))$ for any $x,$ where   $j=2,3.$ If $\st\le \wedge,$ i.e., $a\st b\le a\wedge b$ for any $a,b\in [0,\by],$ then 
\begin{align}\label{b1}
\psi_1(\varphi_1(a\st b)\wedge c\wedge d)\ge \psi_2(\varphi_2(a)\wedge c)\st \psi_3(\varphi_3(b)\wedge d)
\end{align}
for all $a,b\in [0,\by]$ and $c,d\in D.$ 
\end{lemma}

\proof
In a similar manner as in the proof of Theorem~\ref{ctw2} we can show that $\psi$'s are well-defined in \eqref{b1}.
To shorten the notation, we put $$P(a,b,c,d):=\psi_2(\varphi_2(a)\wedge c)\st \psi_3(\varphi_3(b)\wedge d).$$
From the obvious inequalities and assumptions on $\psi_i$ and $\varphi_i$, we have
\begin{align}\label{c7.1}
P(a,b,c,d)&\le \psi_2(\varphi_2(a))\st \psi_3(\varphi_3(b))\le a\st b\le \psi_1(\varphi_1(a\st b))
\end{align}
for all $a,b\in [0,\by]$ and $c,d\in D.$ Let $c,d\in D.$ Consider four cases:
\begin{enumerate}[noitemsep]
\item[(i)] Firstly, suppose $c\le\varphi_2(\by)$ and $d\le\varphi_3(\by).$
Applying $\psi_j\le \psi_1$  and $\st\le \wedge,$ we get
 \begin{align*}
P(a,b,c,d)&\le \psi_2(c)\st \psi_3(d)\le  \psi_1(c\wedge d)
\end{align*}
for any $a,b\in [0,\by].$ By \eqref{c7.1} and monotonicity of $\psi_1$ we obtain \eqref{b1}.
\item[(ii)] Let 
$c>\varphi_2(\by)$ and $d>\varphi_3(\by).$
From \eqref{c7.1}  we conclude that
\begin{align*}
P(a,b,c,d)&\le \psi_1(\varphi_1(a\st b))
=\psi_1(\varphi_1(a\st b) \wedge \varphi_1(\by))
\end{align*}
for any $a,b\in [0,\by].$ By $\varphi_1(\by)=\varphi_j(\by)$ for $j=2,3,$
\begin{align*}
P(a,b,c,d) \le \psi_1(\varphi_1(a\st b) \wedge \varphi_2(\by)\wedge \varphi_3(\by))\le\psi_1(\varphi_1(a\st b) \wedge c\wedge d)
 \end{align*}
 for all $a,b\in[0,\by].$
 \item[(iii)] If $c>\varphi_2(\by)$ and $d \le \varphi_3(\by),$ then
\begin{align}\label{c7.2}
P(a,b,c,d)\le  \psi_2(\varphi_2(a))\wedge \psi_3(d) \le a\wedge \psi_1(d) \le \psi_1(\varphi_1(a)\wedge d)
\end{align}
for any $a,b\in [0,\by]$.
Combining \eqref{c7.1} and \eqref{c7.2}, we get
\begin{align*}
P(a,b,c,d)&\le \psi_1(\varphi_1(a\st b)\wedge \varphi_1(a)\wedge d)
\le  \psi_1 (\varphi_1(a\st b)\wedge \varphi_2(\by)\wedge d)
\\&\le  \psi_1 (\varphi_1(a\st b)\wedge c\wedge d)
\end{align*}
for all $a,b\in [0,\by].$
\item[(iv)] The case $c\le \varphi_2(\by)$ and $d >\varphi_3(\by)$ is similar to (iii).
\end{enumerate}
By (i)--(iv) we get \eqref{b1} for any $a,b\in [0,\by]$ and $c,d\in D.$\qed\medskip

\begin{tw}\label{tw:Sug}
Suppose that $\by\in (0,\infty],$ $m\in\cMm$ and $\st$ is a~non-decreasing and left-continuous fusion function such that $\st \le \wedge.$ Let
$\varphi_i\colon [0,\by]\to [0,\by]$ be functions such that $\varphi_1(\by)=\varphi_j(\by),$ $\varphi_1$ is non-decreasing and $\varphi_j$ are increasing and right-continuous, $\psi_i\colon [0,\varphi_i(\by)]\to [0,\by]$ be non-decreasing and $\psi_j$ be left-continuous such that $\psi_1\ge \psi_j$ and $\psi_j(\varphi_j(x))\le x\le \psi_1(\varphi_1(x))$ for any $x,$  where $i=1,2,3$ and $j=2,3.$ Then  the Chebyshev type inequality 
\begin{align*}
 \psi_1\big(\calka{\rM,A}{ \varphi_1(f\st g)}\big) \ge \psi_2\big(\calka{\rM,A}{ \varphi_2(f)}\big) \st \psi_3\big(\calka{\rM,A}{ \varphi_3(g)}\big)
\end{align*}
is valid for all comonotone functions $f,g\in \cFf$ on $A\in\cA.$
\end{tw}

\proof
Put $k=\by,$ $\loo=\st,$ $B=A,$ $\vtr=\wedge,$  and $\c_i=\wedge$  for $i=1,2,3$ in Theorem~\ref{ctw2}\,(a). From Lemma~\ref{lempom2} with $D=m(\cA)$ the inequality \eqref{b1} is valid for any $a,b\in [0,\by]$ and $c,d\in m(\cA).$  From Theorem~\ref{ctw2}\,(a)  we get the statement.
\qed\medskip

As the~special case with $\psi_i=\varphi_i^{-1}$ and $\varphi_i=\varphi$ we obtain the following result.

\begin{cor}\label{cor:p17.1}
Let $\by\in (0,\infty],$ $m\in\cMm$ 
%with $m(X)\le \by$
and
$\st$ be a~non-decreasing and left-continuous fusion function such that $\st \le \wedge.$ Assume that
$\varphi\colon [0,\by]\to [0,\by]$ is continuous and increasing function such that $\varphi(0)=0,$  where $i=1,2,3$ and $j=2,3.$ Then  the Chebyshev type inequality 
\begin{align}\label{b5}
 \varphi^{-1}&\big(\calka{\rM,A}{ \varphi(f\st g)}\big)\ge \varphi^{-1}\big(\calka{\rM,A}{ \varphi(f)}\big) \st \varphi^{-1}\big(\calka{\rM,A}{ \varphi(g)}\big)
\end{align}
is valid for all comonotone functions $f,g\in \cFf$ on $A\in\cA.$
\end{cor}

%\begin{cor}\label{cor:p17.1}
%Let $\by\in (0,\infty],$ $m\in\cMm$ and $\st$ be a~non-decreasing and left-continuous fusion function such that $\st \le \wedge.$ Let $\varphi_i\colon [0,\by]\to [0,\by]$ be continuous and increasing such that $\varphi_1\le \varphi_j,$ $\varphi_1(\by)=\varphi_j(\by)$ and $\varphi_j(0)=0,$  where $i=1,2,3$ and $j=2,3.$ Then  the Chebyshev type inequality \begin{align}\label{b5}
 %\varphi_1^{-1}&\big(\calka{\rM,A}{ \varphi_1(f\st g)}\big)\ge \varphi_2^{-1}\big(\calka{\rM,A}{ \varphi_2(f)}\big) \st \varphi_3^{-1}\big(\calka{\rM,A}{ \varphi_3(g)}\big)\end{align}
%is valid for all comonotone functions $f,g\in \cFf$ on $A\in\cA.$\end{cor}
The following example shows that the assumption $\varphi(0)=0$ cannot be abandoned from Corollary~\ref{cor:p17.1}.  %\cite[Corollary 4.5]{AEM_Holder} being a consequence of \cite[Theorem 4.1]{AEM_Holder}.

\begin{example}\label{cex1}
Let  $\varphi(x)=0{.}5(x+1)\mI{[0,1]}(x),$ $f=g=0{.}5\mI{X}$  and consider a~capacity $m$ such that $m(A)=0{.}4,$ where $A\in\cA$. Clearly, all assumptions of Corollary~\ref{cor:p17.1} with  $\by=1,$ $\varphi_i=\varphi$ and $\st=\cdot$  are satisfied except $\varphi(0)=0.$ Then
\begin{align*}
\varphi^{-1}\big(\calka{\rM, A} {\varphi(f g)} \big)&\ge \varphi^{-1}\big(\calka{\rM, A} {\varphi(f)}\big)\cdot \varphi^{-1}\big(\calka{\rM, A}{\varphi(g)}\big),\\
\varphi^{-1}(0{.}625\wedge 0{.}4)&\ge \varphi^{-1}(0{.}75\wedge 0{.}4)\cdot \varphi^{-1}(0{.}75\wedge 0{.}4).
\end{align*}
However, the value $\varphi^{-1}(0{.}4)$ is not defined, so the Chebyshev type inequality
%\cite[Corollary 4.5]{AEM_Holder} 
cannot hold. 
% The same corollary can be found in \cite[Corollary 3.9]{AMOPS_Chebyshev} or \cite[Corollary 3.20]{AMV_predictive}. 
\end{example}

\begin{remark}\label{rem:false}
In \cite[Theorem 3.1]{AgahiMesiarOuyang09} 
%and \cite[Theorem 4.1]{AgahiMesiarOuyang12}
the Chebyshev type inequality of the form \eqref{b5}  for a~fixed $A\in\cA$  has been obtained without the assumption on $\varphi(0)=0.$ 
 %TU będą argumenty na kontratak
 %However, Theorem 3.1 from \cite{AgahiMesiarOuyang09} will be true if we add an additional assumption  $\varphi(0)\le m(A).$ To see it,  use the method of the proof presented in \cite{AgahiMesiarOuyang09} and observe that $\calka{\rM,A}{\varphi(h)}\in [\varphi(0),\sup_x \varphi(x)]$ for any $h\in\cFf[\infty],$ since $\calka{\rM,A}{\varphi(h)}\ge \varphi(0)\wedge m(A)=\varphi(0).$
\end{remark}
  
As a~consequence of the Chebyshev type inequality, we get the following Liapunov type inequality. 

\begin{cor}
 Let $\by\in (0,\infty],$ $m\in\cMm$  and $\varphi_i\colon [0,\by]\to [0,\by]$ be increasing and right-continuous, $\psi_i\colon [0,\varphi_i(\by)]\to [0,\by],$ $i=1,2,$ such that  $\varphi_1(\by)=\varphi_2(\by),$ $\psi_1\ge \psi_2$ and $\psi_2(\varphi_2(x))\le x\le \psi_1(\varphi_1(x))$ for any $x.$ Then the Liapunov type inequality 
\begin{align*}
 \psi_1&\big(\calka{\rM,A}{ \varphi_1(f)}\big)\ge \psi_2\big(\calka{\rM,A}{ \varphi_2(f)}\big)
\end{align*}
holds for each $f\in \cFf$ and $A\in\cA.$
\end{cor}
\proof
Put $\st=\wedge,$  $g=f,$ $\varphi_3=\varphi_2,$ 
and $\psi_3=\psi_2$ in Theorem~\ref{tw:Sug}.
\qed\medskip

\section{The Chebyshev type inequality for any functions}\label{sub:further}

In this section, we show that, due to the definition of $m$-positively dependent functions, one can derive the Chebyshev type inequality for any functions strengthening assumptions on monotone measure.

\begin{tw}\label{tw:51}
Assume that %$\by>0,$ 
$\by\in(0,\infty],$ $k\in (0,\by]$ and  $m\in\cMm$ such that $m(C\cap D)\ge m(C)\vtr m(D)$ for all $C,D\in\cA$ with $\vtr\colon m(\cA)\times m(\cA)\to m(\cA).$  Let
$\varphi_i\colon [0,\by]\to [0,\by]$ be increasing and right-continuous, and $\psi_i\colon [0,\varphi_i(\by)]\to [0,\by]$ be non-decreasing and left-continuous for $i=1,2,3.$
Suppose that  $\st$ is a~non-decreasing and left-continuous fusion function, $\c_i\colon [0,\by]\times m(\cA)\to [0,\by]$ are non-decreasing such that $\varphi_i(\by)\c_i m(X)\le \varphi_i(\by)$ and $\by\c_j 0=0$ for $i=1,2,3$ and $j=2,3.$
If \eqref{c2a} with $\loo=\st$ holds for any $a,b\in [0,k]$ and $c,d\in m(\cA)$,
then the Chebyshev type inequality \eqref{c2b}
is satisfied for all functions $f|_A,g|_B\in \cFf[k]$.
\end{tw}

\proof
Observe that all functions $f|_A$ and $g|_B$ are $m$-positively dependent with respect to $\vtr,$ since $m(C\cap D)\ge m(C)\vtr m(D)$ for all $C,D\in\cA.$ Applying
Theorem~\ref{ctw2}\,(a) with $\loo=\st$ we get the statement.
\qed \medskip

Using Theorem~\ref{tw:51}, we can obtain many versions of Chebyshev type inequality for any functions. Below, we present one proposition. 
%We say that $m\in\cMm$ is \textit{minitive} if $m(A\cap B)=m(A)\wedge m(B)$ for any $A,B\in\cA.$
%When replacing $\c_i$ with $\wedge$ in Theorem~\ref{tw:51}, we obtain the Chebyshev type inequality for the Sugeno integral. 
%But, unfortunately, the condition \eqref{c2a} with $\loo=\st$ and $\rS=\rM$ is very hard to check in many situations. But thanks to Lemma~\ref{lempom2} we may change it to another condition which is much more easier to verify.

\begin{cor}\label{cor:11}
Let $m\in\cMm$ be minitive and $\st$ be a~non-decreasing left-continuous fusion function such that $\st\le \wedge.$ Let
$\varphi_i\colon [0,\by]\to [0,\by]$ be increasing and right-continuous and
$\psi_i\colon [0,\varphi_i(\by)]\to [0,\by]$ be non-decreasing and left-continuous, where $i=1,2,3.$ Assume that $\varphi_1(\by)=\varphi_j(\by)$, $\psi_1 \ge \psi_j$  and $\psi_j(\varphi_j(x))\le x\le \psi_1(\varphi_1(x))$ for all $x$ and $j=2,3.$
Then the Chebyshev type inequality
\begin{align*}
\psi_1\big(\calka{\rM,A\cap B}{\varphi_1(f\st g)}\big)\ge \psi_2\big(\calka{\rM,A}{\varphi_2(f)}\big)\st \psi_3\big(\calka{\rM,B}{\varphi_3(g)}\big)
\end{align*}
is fulfilled for all functions $f|_A,g|_B\in \cFf$.
\end{cor}
\proof
Put $\c_i=\vtr=\wedge$ and $k=\by$ in Theorem~\ref{tw:51}. By assumptions we obtain the inequality \eqref{c2a} for all $a,b\in [0,\by]$ and $c,d\in m(\cA)$  from Lemma~\ref{lempom2} with $D=m(\cA)$ and $\c_i=\wedge.$ Then Theorem  \ref{tw:51} provides the statement.
\qed\medskip

%Moreover, we can state the analogous result for the seminormed fuzzy integral.

The following examples illustrate the above results. 

\begin{example}
Let $X=[0,1]$ and consider on $(X,\cA)$ the minitive capacity
%\footnote{Monotone measure $m$ is minitive, if $m(C\cap D)=m(C)\wedge m(D)$ for any measruables sets $C,D.$}
$m(A)=1-\sup_{x\in A^c} x.$ Assume that $\st=\cdot,$ $\psi_i(x)=\sqrt{x}$, $\varphi_j(x)=x^2,$ $\varphi_1(x)=x,$ $f(x)=\mI{[0,\,0.5]}(x)+0.5\mI{(0.5,\,1]}(x)$ and $g(x)=x$, where $i=1,2,3$ and $j=2,3.$ 
Fix $A=B=X$. Then all assumptions in Corollary~\ref{cor:11} are valid with $\by=1$. By Corollary~\ref{cor:11} the Chebyshev type inequality has the form
\begin{align*}
\calka{\rM}{x\cdot\mI{[0,\,0.5]}(x)+0.5\,x\cdot\mI{(0.5,\,1]}(x)}\ge \calka{\rM}{\mI{[0,\,0.5]}(x)+0.25\cdot\mI{(0.5,\,1]}(x)}\cdot \calka{\rM}{x^2}.
\end{align*}
By simple calculation, we get
\begin{align*}
\calka{\rM}{x\,\mI{[0,\,0.5]}(x)+0.5\,x\cdot\mI{(0.5,\,1]}(x)}&=\sup_{t\in [0,\,0.25]}\{t\wedge m([t,1])\}\vee \sup_{t\in(0.25,\,0.5]}\{t\wedge m([t,0.5]\cup [2t,1])\}
\\&=\sup_{t\in [0,\,0.25]}\{t\wedge (1-t)\}\vee \sup_{t\in(0.25,\,0.5]}\{t\wedge (1-2t)\}=1/3.
\end{align*}
In a~similar manner we obtain
\begin{align*}
\calka{\rM}{\mI{[0,\,0.5]}(x)+0.25\cdot\mI{(0.5,\,1]}(x)}&=\sup_{t\in [0,\,0.25]} \{t\wedge m(X)\}\vee \sup_{t\in(0.25,\,1]}\{t\wedge 0\}=0.25,\\
\calka{\rM}{x^2}&=\sup_{t\in [0,1]} \{t\wedge (1-\sqrt{t})\}=(3-\sqrt{5})/2.
\end{align*}
%If $\varphi_j(x)=\sqrt{x}$ for $j=2,3,$ then the assumptions in Corollary~\ref{cor:11} are not satisfied. Additionally, the Chebyshev type inequality does not hold, since
%\begin{align*}
%\calka{\rM}{\sqrt{f}}&=\sup_{t\in [0,1]} \{t\wedge (1-t^2)\}=(\sqrt{5}-1)/2,\\
%\calka{\rM}{\sqrt{g}}&=\sup_{t\in [0,\sqrt{0.5}]} \{t\wedge m(X)\}\vee \sup_{t\in(\sqrt{0.5},1]}\{t\wedge 0\}=\sqrt{0.5}
%\end{align*}and $0.4370\approx \calka{\rM}{\sqrt{f}}\cdot \calka{\rM}{\sqrt{g}} >1/3.$
\end{example}

Let $\rA,\rB\colon [0,1]^2\to[0,1]$ be two fusion functions. We say that 
\textit{$\rB$ is dominated by $\rA$}, 
if the inequality $\rA(\rB(a,b), \rB(c,d))\ge \rB(\rA(a,c),\rA(b,d))$ holds for all $a,b,c,d \in [0,1].$

\begin{example}\label{ex:52}
Let $m$ is the Lebesgue measure on $X=[0,1].$ Consider $\st =\vtr=\rW,$ $k=\by=1,$ $\c_i=\wedge,$ $\varphi_i(x)=x,$  $\psi_1(x)=\sqrt{x}$ and $\psi_j(x)=x^2$ for $i=1,2,3$ and $j=2,3$ in Theorem~\ref{tw:51}. 
Note that $m(C\cap D)\ge \rW(m(C),m(D))$ for any $C,D\in\mathcal{A}$ since $m$ is supermodular, see Example~\ref{ex:rec1}.
The inequality \eqref{c2a} with $\loo=\st$ holds by the reason of $\rW$ being dominated by $\rM.$ 
By Theorem~\ref{tw:51} with $A=B=X$ the Chebyshev type inequality 
\begin{align}\label{b10}
    \sqrt{\calka{\rM}{\rW(f,g)}}\ge \rW\big(\big(\calka{\rM}{f}\big)^2, \big(\calka{\rM}{g}\big)^2\big)
\end{align}
is true for any $f,g\in \cFf[1].$ For $f(x)=2x^2-2x+1$ and $g(x)=-2x^2+2x$ we get $\calka{\rM}{\rW(f,g)}=0$ and
\begin{align*}
\calka{\rM}{f}&=0.5\vee\sup_{t\in (0.5,\,1]} \{t\wedge (1-2\sqrt{0.5(t-0.5)})\}=2-\sqrt{2},\\
\calka{\rM}{g}&=\sup_{t\in (0,\,0.5]} \{t\wedge 2\sqrt{0.5(0.5-t)}\}=\sqrt{2}-1,
\end{align*}
so the inequality in \eqref{b10} becomes equality.
\end{example}
\noindent Example~\ref{ex:52} illustrates that the equality in the Chebyshev type inequality may be achieved by non-constant
functions $f$ and $g$ (cf. Examples \ref{ex:rown1}--\ref{ex:rown2}).

\section*{Conclusion} 
In the present paper, we have focused on Chebyshev type inequalities for generalized (upper) Sugeno integral for $m$-positively dependent functions, and for any comonotone functions and monotone measure. We have presented a~technique how to obtain some known results from the literature. Moreover, we have provided a~few unknown results in the literature so far. For instance, when restricting the considered fusion functions to fuzzy conjunctions, we can state general Chebyshev type inequalities for q-integral recently introduced and studied by Dubois et al. in~\cite{DuboisPradeRicoTeheux17}.
%see also the recent paper~\cite{YanOuyang19} on the Chebyshev type inequality for q-integral. 
In the point of view that the considered integrals are aggregation functions, we expect applications of our results everywhere where some bounds on the aggregation process is needed, such as information aggregation, or decision making. 

\section*{Acknowledgement}
This work was supported by the Slovak Research and Development Agency
under the  contract  No.  APVV-16-0337. The work is also cofinanced by bilateral call Slovak-Poland
grant scheme No. SK-PL-18-0032 together with
the Polish National Agency for Academic Exchange  PPN/BIL/2018/1/00049/U/00001.

\end{document}